\documentclass[11pt]{article}

\usepackage{morefloats}
\usepackage{amsfonts}
\usepackage{color}
\usepackage{multirow}
\usepackage{latexsym}
\usepackage{amsmath,amssymb,amsthm,graphicx}
\usepackage{dsfont}
\usepackage{extarrows}
\usepackage{hyperref}

\numberwithin{equation}{section}

\newtheoremstyle{thm}
{9pt}
{9pt}
{\itshape}
{}
{\bfseries}
{.}
{ }
{}
\theoremstyle{thm}

\newtheorem{theorem}{Theorem}[section]

\newtheoremstyle{def}
{9pt}
{9pt}
{}
{}
{\bfseries}
{.}
{ }
{}
\theoremstyle{def}

\newcommand{\PP}{\mathbb{P}}

\newcommand{\stktheta}{\stackrel{\mbox{\scriptsize$\PP_\vartheta$}}{\longrightarrow}}
\newcommand{\R}{\mathbb{R}} 
\newcommand{\N}{\mathbb{N}} 
\newcommand{\E}{\mathbb{E}} 

\renewcommand{\footnoterule}{%
	\kern -3.5pt
	\hrule width \textwidth height 1pt
	\kern 3.5pt
}

\makeatletter
\def\blfootnote{\xdef\@thefnmark{}\@footnotetext}
\makeatother

\begin{document}

\title{\bf  Bahadur efficiencies of the Epps--Pulley test for normality}


\author{Bruno Ebner and Norbert Henze}

\date{\today}
\maketitle

\begin{center}
{\Large In memoriam Yakov Yu. Nikitin (1947--2020)}
\end{center}

\blfootnote{
{\em Test for normality; empirical characteristic function; Bahadur efficiency; eigenvalues of integral operators} }

\begin{abstract}
The test for normality suggested  by Epps and Pulley \cite{EP:1983} is a serious competitor to tests based on the empirical distribution function.
In contrast to the latter procedures, it has been generalized to obtain a genuine affine invariant and universally consistent test for normality
in any dimension. We obtain approximate Bahadur efficiencies for the test of Epps and Pulley, thus complementing recent results of Milo\u{s}evi\'{c} et al.
(see \cite{MNO:2021}). For certain values of a tuning parameter that is inherent in the Epps--Pulley test, this test outperforms
 each of its  competitors considered in \cite{MNO:2021}, over the whole range of six close alternatives to normality.
\end{abstract}

\noindent
\section{Introduction}\label{sec:Intro}
The purpose of this article is to derive Bahadur efficiencies for
the test of normality proposed by Epps and Pulley \cite{EP:1983}, thus complementing recent results of
Milo\u{s}evi\'{c} et al.  \cite{MNO:2021}, who confined their study to tests of normality based on the empirical distribution
function. To be specific, suppose $X_1,X_2, \ldots$ is a sequence of
independent and identically distributed (i.i.d.) copies of a random variable $X$ that has an absolutely continuous
distribution with respect to Lebesgue measure. To test the hypothesis $H_0$ that the distribution  of $X$ is
some unspecified non-degenerate normal distribution, Epps and Pulley \cite{EP:1983} proposed to use the test statistic
\[
T_{n,\beta} = n \int_{-\infty}^\infty \Big{|} \psi_n(t) - {\rm e}^{-t^2/2}\Big{|}^2 \, \varphi_\beta(t)\, {\rm d}t.
\]
Here, $\psi_n(t) = n^{-1}\sum_{j=1}^n \exp\big({\rm i}tY_{n,j}\big)$ is the empirical characteristic function of
the so-called {\em scaled residuals} $Y_{n,1}, \ldots, Y_{n,n}$, where $Y_{n,j} = S_n^{-1} (X_j - \overline{X}_n)$, $j=1,\ldots,n$,
and $\overline{X}_n = n^{-1} \sum_{j=1}^n X_j$, $S_n^2 = n^{-1}\sum_{j=1}^n (X_j- \overline{X}_n)^2$ are the sample mean and
the sample variance of $X_1,\ldots,X_n$, respectively, and $\beta>0$ is a so-called \textit{tuning parameter}. Moreover,
\[
\varphi_\beta(t) = \frac{1}{\beta \sqrt{2\pi}} \exp \Big( - \frac{t^2}{2\beta^2 } \Big), \quad t \in \mathbb{R},
\]
is the density of the centred normal distribution with variance $\beta^2$. A closed-form expression of $T_{n,\beta}$ that is amenable to
computational purposes is
\begin{equation}\label{deftnb}
T_{n,\beta} = \frac{1}{n} \sum_{j,k=1}^n \exp \! \bigg( \! \! - \frac{\beta^2}{2}\big(Y_{n,j}\! - \! Y_{n,k}\big)^2 \! \bigg) \!
 - \frac{2}{\sqrt{1 \! + \! \beta^2}} \sum_{j=1}^n \exp \! \bigg( \! \! - \frac{\beta^2 Y_{n,j}^2 }{2(1\! + \! \beta^2)} \! \bigg)
 + \frac{n}{\sqrt{1\! + \! 2\beta^2}}.
\end{equation}

Epps and Pulley did not obtain neither the limit null distribution of $T_{n,\beta}$ as $n \to \infty$
nor the consistency of a test for normality that rejects $H_0$ for large values of $T_{n,\beta}$. Their procedure,
however, turned out to be a serious competitor to the classical tests of Shapiro--Wilk, Shapiro--Francia and
Anderson--Darling in simulation studies (see \cite{BDH:1989}). In the special case $\beta =1$, Baringhaus and Henze \cite{BH:1988} generalized
the approach of Epps and Pulley to obtain a genuine test of multivariate normality, and they derived the limit null distribution
of $T_{n,1}$. Moreover, they proved the consistency of the test of Epps and Pulley against each alternative to normality having a finite second moment.
The latter restriction was removed by S. Cs\"org\H{o} \cite{CS:1989}.
By an approach different from that adopted in \cite{BH:1988}, Henze and Wagner \cite{HW:1997} obtained both the limit null distribution and the
limit distribution of $T_{n,\beta}$ under contiguous alternatives to normality. Under fixed alternatives to normality, the limit distribution
of $T_{n,\beta}$ is normal, as elaborated by \cite{BEH:2017} in much greater generality for weighted $L^2$-statistics. For more information
on $T_{n,\beta}$, especially on the role of the tuning parameter $\beta$, see Section 2.2 of \cite{EH:2020}.

Notice that $T_{n,\beta}$ is invariant with respect to affine transformations
$X_j \mapsto aX_j+b$, where $a,b \in \R$ and $a \neq 0$. Hence, under $H_0$, both the finite-sample and the asymptotic distribution of $T_{n,\beta}$ do not depend
on the parameters $\mu$ and $\sigma^2$ of the underlyling normal distribution N$(\mu,\sigma^2)$. Under $H_0$, we will  thus assume
$\mu=0$ and $\sigma^2=1$. The rest of the paper unfolds as follows: In Section \ref{sec:ABE},  we revisit the notion of approximate Bahadur efficiency.
Sections \ref{sec:stochlimit} and \ref{sec:local} deal with stochastic limits and local Bahadur slopes,  and Section \ref{sec:approxev} tackles an
eigenvalue problem connected with the limit null distribution of the test statistic.  The final section \ref{sec:alternatives} contains
results regarding local approximate Bahadur efficiencies of the Epps--Pulley test for the six close alternatives considered in
\cite{MNO:2021} and a wide spectrum of values of the tuning parameter $\beta$.

\section{Approximate Bahadur efficiency}\label{sec:ABE}
There are several options to compare different tests for the same testing problem as the sample size $n$ tends to infinity, see \cite{NI:1995}.
One of these options is asymptotic efficiency due to Bahadur (see \cite{BA:1967}). This notion of asymptotic efficiency requires knowledge of the
large deviation function of the test statistic. Apart from the notable exception given in \cite{NT:2007}, such knowledge, however, is  hitherto
not available for statistics that contain estimated parameters, like $T_{n,\beta}$ given in \eqref{deftnb}. To circumvent this drawback,
one usually employs the so-called {\em approximate} Bahadur efficiency, which only requires results on the tail behavior of the limit distribution of the
test statistic under the null hypothesis.  To be more specific with respect to the title of this paper, let $X,X_1,X_2, \ldots$ be a sequence of  i.i.d. random variables, where
the distribution of $X$ depends on a real-valued parameter  $\vartheta \in \Theta$, where $\Theta$ denotes the parameter space,
and only the case $\vartheta =0$ corresponds to the case that the distribution
of $X$ is standard normal. Suppose, $(S_n)_{n\ge 1}$, where $S_n = S_n(X_1,\ldots,X_n)$, is a sequence of test statistics of
the hypothesis $H_0: \vartheta =0$ against the alternative $H_1$: $\vartheta \in \Theta \setminus \{0\}$. Furthermore, suppose that rejection of $H_0$
is for large values of $S_n$. The sequence $(S_n)$ is called a {\em standard sequence}, if the following conditions hold (see, e.g., \cite{NI:1995}, p.10, or
\cite{EHN:2009}, p. 3427):
\begin{itemize}
\item
There is a continuous distribution function $G$ such that, for $\vartheta =0$,
\begin{equation}\label{Bcondition1}
\lim_{n\to \infty} \PP_0(S_n \le x) = G(x), \quad x \in \mathbb{R}.
\end{equation}
\item There is a constant $a_S$, $ 0 < a_S <  \infty$, such  that
\begin{equation}\label{Bcondition2}
\log(1-G(x)) = -\frac{a_Sx^2}{2}(1+o(1)) \ \text{as} \ x \to \infty.
\end{equation}
\item There is a real-valued function $b_S(\vartheta)$ on $\Theta \setminus\{0\}$, with $ 0 < b_S(\vartheta) < \infty$, such that, for each $\vartheta \in \Theta \setminus \{0\}$,
\begin{equation}\label{Bcondition3}
\frac{S_n}{\sqrt{n}} \stktheta b_S(\vartheta).
\end{equation}
\end{itemize}
Then the so-called {\em approximate Bahadur slope}
\[
c_S^*(\vartheta) = a_S \cdot b_S^2(\vartheta), \quad \vartheta \in \Theta \setminus \{0\},
\]
is a measure of approximate Bahadur efficiency.  Usually, it is true that
$
c_S^*(\vartheta) \sim \ell(S) \cdot \vartheta^2
$
as $\vartheta \to 0$. In this case $\ell(S)$ is called the {\em local (approximate) index}  of the sequence $(S_n)$. We will see that the sequence
$(S_n)$, where $S_n := \sqrt{T_{n,\beta}}$, is a standard sequence. To this end, we will derive the stochastic limit of $T_{n,\beta}/n$
for a general alternative in Section \ref{sec:stochlimit}. In Section \ref{sec:local}, we will specialize this stochastic limit for local alternatives,
and we will derive the local index for the  Epps--Pulley test statistic.

\section{Stochastic limit of $T_{n,\beta}/n$}\label{sec:stochlimit}
To calculate the asymptotic Bahadur efficiency of the test of Epps and Pulley, we need the following result.

\begin{theorem}\label{thmstk}
Suppose that $\mathbb{E}(X^2) < \infty$. Then
\[
\frac{T_{n\beta}}{n} \stackrel{\mbox{\scriptsize$\PP$}}{\longrightarrow}  \mathbb{E}\bigg[ \exp\left(- \frac{\beta^2(Y_1-Y_2)^2}{2}\right)\bigg] - \frac{2}{\sqrt{1+\beta^2}}\mathbb{E}\bigg[ \exp\left(
- \frac{\beta^2Y_1^2}{2(1\! + \! \beta^2)} \right)\bigg]  + \frac{1}{\sqrt{1\! + \! 2\beta^2}}.
\]
Here, $\stackrel{\mbox{\scriptsize$\PP$}}{\longrightarrow}$ denotes convergence in probability, and $Y_j = (X_j-\mu)/\sigma$, $j \ge 1$, where $\mu = \mathbb{E}(X)$ and $\sigma^2 = \mathbb{V}(X)$.
\end{theorem}

{\sc Proof.} From \eqref{deftnb}, we have
\begin{eqnarray*}
\frac{T_{n,\beta}}{n} & = & \frac{1}{n^2} \sum_{j,k=1}^n \exp \! \bigg( \! \! - \frac{\beta^2}{2}\left(\frac{X_j-X_k}{S_n}\right)^2 \! \bigg) \\
 & & - \frac{2}{\sqrt{1 \! + \! \beta^2}} \cdot \frac{1}{n} \sum_{j=1}^n \exp \! \bigg( \!  - \frac{\beta^2}{2(1\! + \! \beta^2)} \left(\frac{X_j- \overline{X}_n}{S_n} \right)^2 \! \bigg)
 + \frac{1}{\sqrt{1\! + \! 2\beta^2}}\\
& =: & A_{n,1}  - \frac{2}{\sqrt{1 \! + \! \beta^2}} \cdot A_{n,2} +  \frac{1}{\sqrt{1\! + \! 2\beta^2}}
\end{eqnarray*}
(say). By symmetry, it follows that
\begin{eqnarray*}
\mathbb{E}\big(A_{n,1}\big) & = & \frac{1}{n} + \frac{n-1}{n} \cdot \mathbb{E}\bigg[ \exp \! \bigg( \! \! - \frac{\beta^2}{2}\left(\frac{X_1-X_2}{S_n}\right)^2 \! \bigg],\\
\mathbb{E}\big(A_{n,2}\big) & = &  \mathbb{E} \bigg[\exp \! \bigg( \!  - \frac{\beta^2}{2(1\! + \! \beta^2)} \left(\frac{X_1- \overline{X}_n}{S_n} \right)^2 \! \bigg) \bigg].
\end{eqnarray*}
Since $\overline{X}_n \rightarrow \mu$ and $S_n \rightarrow \sigma$ almost surely as $n \to \infty$ by the strong law of large numbers, it follows from
Lebesgue's dominated convergence theorem that
\[
\lim_{n \to \infty} \mathbb{E}\big(A_{n,1}\big) = \mathbb{E}\bigg[ \exp\left(- \frac{\beta^2(Y_1-Y_2)^2}{2}\right)\bigg], \quad
\lim_{n \to \infty} \mathbb{E}\big(A_{n,2}\big) = \mathbb{E}\bigg[ \exp\left(
- \frac{\beta^2Y_1^2}{2(1\! + \! \beta^2)} \right)\bigg],
\]
and thus the expectation of $T_{n,\beta}/n$ converges to the stochastic limit figuring in Theorem \ref{thmstk}. Likewise, the variance
of  $T_{n,\beta}/n$ converges to zero.

\bigskip

\section{Local Bahadur slopes}\label{sec:local}
As was done in Milo\u{s}evi\'{c} et al. \cite{MNO:2021}, we now assume that ${\cal G} = \{G(x;\vartheta)\}$ is a family of distribution functions (DF's)
with densities $g(x;\vartheta)$, such that $\vartheta =0$ corresponds to the standard normal DF $\Phi$ and density $\varphi$, and each of the distributions
for $\vartheta \neq 0$ is non-normal. Moreover, we  assume that the regularity assumptions WD in \cite{NP:2004} are satisfied.
If $X,X_1, X_2, \ldots $ are i.i.d. random variables with DF $G(\cdot;\vartheta)$, we have to consider the stochastic limit figuring in Theorem \ref{thmstk}
as a function of $\vartheta$ and expand this function at $\vartheta =0$. To this end, let
\begin{equation}\label{defdeltagam}
\gamma = \frac{\beta^2}{2}, \ \delta = \frac{\beta^2}{2(1+\beta^2)}.
\end{equation}
Then, putting
\begin{eqnarray*}
\mu(\vartheta) & = & \int x g(x;\vartheta) \, \text{d}x,\\
\sigma^2(\vartheta) & = & \int x^2  g(x;\vartheta) \, \text{d}x - \mu^2(\vartheta),
\end{eqnarray*}
Theorem \ref{thmstk} yields
\[
\frac{T_{n,\beta}}{n} \stktheta b_{T_\beta}(\vartheta),
\]
where $\stktheta$ denotes convergence in probability under the true parameter $\vartheta$, and
\begin{eqnarray}\label{integ1}
b_{T_\beta}(\vartheta) & = & \iint \exp\left(-  \frac{\gamma(x-y)^2}{\sigma^2(\vartheta)} \right) g(x;\vartheta) g(y;\vartheta) \, \text{d}x \, \text{d} y\\ \label{integ2}
& & \hspace*{5mm} - \frac{2}{\sqrt{1+\beta^2}} \, \int \exp\left(- \frac{\delta (x-\mu(\vartheta))^2}{\sigma^2(\vartheta)} \right) g(x;\vartheta)\, \text{d} x + \frac{1}{\sqrt{1+2\beta^2}}.
\end{eqnarray}
Here and in what follows, each unspecified integral is over $\R$.

Notice that $b_{T_\beta}(0) =0$. We have to find the quadratic (first non-vanishing) term in the Taylor expansion of  $b_{T_\beta}$ around zero, i.e, we look for some (local index) $\Delta_\beta >0$ such that
\[
b_{T_\beta}(\vartheta) = \Delta_\beta \, \vartheta^2  + o(\vartheta^2) \quad \text{as} \ \vartheta \to 0.
\]
Writing $g'_\vartheta (x;\vartheta)$, $g''_\vartheta (x;\vartheta)$  for derivatives of $g(x;\vartheta)$ with respect to $\vartheta$, we have
\[
g(x;\vartheta) = \varphi(x) + \vartheta \cdot g'_\vartheta(x;0) + \frac{\vartheta^2}{2}g''_\vartheta(x;0) + O(\vartheta^3)
\]
and thus -- since $\mu(0) =0$ and $\sigma^2(0) = 1$ --
\begin{eqnarray*}
\mu(\vartheta) & = & \vartheta \int xg'_\vartheta(x;0)\, \text{d}x + \frac{\vartheta^2}{2} \int xg''_\vartheta(x;0)\, \text{d} x  + O(\vartheta^3),\\
\sigma^2(\vartheta) & = & 1 + \vartheta \int x^2g'_\vartheta(x;0)\, \text{d}x + \frac{\vartheta^2}{2} \int x^2g''_\vartheta(x;0)\, \text{d} x - \mu(\vartheta)^2  + O(\vartheta^3).
\end{eqnarray*}
Consequently, putting
\begin{equation}\label{abkmusig}
\mu_1 := \mu'(0), \  \mu_2:= \mu''(0), \ \sigma_1 := (\sigma^2)'(0), \ \sigma_2 := (\sigma^2)''(0)
\end{equation}
for the sake of brevity, it follows that
\begin{eqnarray*}
\mu_1 & = & \int x g'_\vartheta(x;0)\, \text{d}x, \quad \mu_2 = \int x g''_\vartheta(x;0) \, \text{d}x, \\
\sigma_1 & = & \int x^2 g'_\vartheta(x;0)\, \text{d}x, \quad \sigma_2 = \int x^2g''_\vartheta(x;0) \, \text{d}x  - 2 \mu'(0)^2.
\end{eqnarray*}
To tackle the integral that figures in \eqref{integ1}, notice that
\begin{eqnarray*}
g(x;\vartheta)g(y;\vartheta) & = & \varphi(x)\varphi(y) + \vartheta\big[g'_\vartheta(x;0)\varphi(y) + g'_\vartheta(y;0)\varphi(x)\big]\\
& & + \vartheta^2 \Big[\frac{1}{2}g''_\vartheta(x;0)\varphi(y) + \frac{1}{2}g''_\vartheta(y;0)\varphi(x) + g'_\vartheta(x;0)g'_\vartheta(y;0)\Big] + O(\vartheta^3).
\end{eqnarray*}
Moreover, it follows from a geometric series expansion that
\begin{equation}\label{exp1s}
\frac{1}{\sigma^2(\vartheta)}  =  1 - \vartheta \sigma_1 + \vartheta^2 \Big[ \sigma_1^2 - \frac{\sigma_2}{2} \Big] + O(\vartheta^3)
\end{equation}
(say). From an expansion of the exponential function, we thus obtain
\begin{eqnarray*}
\! & \! \!  & \! \exp\left(- \frac{\gamma (x-y)^2}{\sigma^2(\vartheta)}\right)  \\
\! & \! = \! & \! {\rm e}^{-\gamma (x-y)^2} \Big[ 1 + \vartheta \sigma_1  \gamma (x\! -\! y)^2 + \vartheta^2 \Big\{ \frac{1}{2}\, \sigma_1^2 \gamma^2 (x\! -\! y)^4 -
\Big(\sigma_1^2 \! - \! \frac{\sigma_2}{2} \Big) \gamma (x\! -\! y)^2 \Big\}   \Big]+ O(\vartheta^3).
\end{eqnarray*}
Using
\[
\iint {\rm e}^{-\gamma (x-y)^2} (x-y)^{2k} \varphi(x)\varphi(y)\, \text{d}x \text{d} y = \frac{4^k \Gamma(k+1/2)}{\sqrt{\pi}(4\gamma +1)^{k+1/2}}, \quad k=0,1,2,
\]
\[
\int {\rm e}^{-\gamma(x-y)^2} \varphi(x)\, \text{d}x = \frac{1}{\sqrt{1+\beta^2}} \cdot {\rm e}^{-\delta y^2}
\]
\[
\int {\rm e}^{-\gamma (x-y)^2} (x-y)^2 \varphi(y)\, \text{d}y = {\rm e}^{-\delta x^2} \cdot \frac{x^2+\beta^2 + 1}{(1+\beta^2)^{5/2}},
\]
and putting
\begin{equation}\label{defdj}
D_0 = \iint {\rm e}^{-\gamma(x-y)^2} g'_\vartheta(x;0)g'_\vartheta(y;0)\, \text{d}x\text{d}y,
\end{equation}
\begin{equation}\label{defjk}
J_{1,k}  =  \int \! {\rm e}^{-\delta x^2} x^k g'_\vartheta(x;0) \, \text{d} x, \quad k = 0,1,2; \qquad J_2   =  \int \! {\rm e}^{-\delta x^2} g''_\vartheta(x;0) \, \text{d} x,
\end{equation}
some algebra gives
\begin{eqnarray}\nonumber
& & \iint \exp\left(-  \frac{\gamma(x-y)^2}{\sigma^2(\vartheta)} \right) g(x;\vartheta) g(y;\vartheta) \, \text{d}x \, \text{d} y \\ \label{erstesint}
& = & \frac{1}{\sqrt{1+2\beta^2}} + \vartheta \bigg\{ \frac{2J_{1,0}}{\sqrt{1+\beta^2}}  + \frac{2 \sigma_1  \gamma}{(1+2\beta^2)^{3/2}}\bigg\} \\ \nonumber
& & + \vartheta^2 \bigg\{ \frac{J_2}{\sqrt{1+\beta^2}}  + D_0  + \frac{2\sigma_1 \gamma \big(J_{1,2} + (\beta^2\! + \! 1) J_{1,0}\big)}{(1+\beta^2)^{5/2}}
+ \frac{\beta^2\big((2\beta^2+1)\sigma_2-(\beta^2+2)\sigma_1^2 \big)}{2(1+2\beta^2)^{5/2}} \bigg\}.
\end{eqnarray}

As for the integral figuring in \eqref{integ2}, we use \eqref{exp1s}. Neglecting terms  that are of order $O(\vartheta^3)$, straightforward but tedious calculations give
\[
\exp\left(\! -\frac{\delta(x\! - \! \mu(\vartheta))^2}{\sigma^2(\vartheta)}\right) =  {\rm e}^{-\delta x^2}\cdot\bigg\{ 1+\delta \vartheta U(x) + \frac{\delta \vartheta^2}{2} V(x) \bigg\} + O(\vartheta^3),
\]
where -- recalling \eqref{abkmusig}  --
\begin{eqnarray*}
U(x) & = &  \sigma_1x^2+2\mu_1x,\\
V(x) & = & \delta \sigma_1^2x^4 \! + \! 4\delta\mu_1\sigma_1 x^3 \! + \! \big(4\delta \mu_1^2-2\sigma_1^2 \! +\! \sigma_2 \big) x^2  \! - \!
(4\mu_1\sigma_1 \! -\! 2\mu_2)x \! - \! 2 \mu_1^2.
\end{eqnarray*}

Thus,
\begin{eqnarray*}
 \int \exp\left(- \frac{\delta (x-\mu(\vartheta))^2}{\sigma^2(\vartheta)}\right) g(x;\vartheta) \, \text{d} x
\! \! & \! \! = \! \! & \! \! \int \! {\rm e}^{-\delta x^2}\Big(1 \! + \!  \delta \vartheta U(x) \! + \!  \frac{\delta \vartheta^2}{2} V(x) \Big) \varphi(x) \, \text{d} x \\
\! \!  &\! \! \! \!  & \! \! + \vartheta \int \! {\rm e}^{-\delta x^2}\Big(1 \! + \! \delta \vartheta U(x) \! + \! \frac{\delta \vartheta^2}{2} V(x) \Big)  g'_\vartheta(x;0) \, \text{d} x \\
\! \!  &\! \! \! \!  & \! \! + \frac{\vartheta^2}{2} \int \! {\rm e}^{-\delta x^2}\Big(1 \! + \! \delta \vartheta U(x) \! + \! \frac{\delta \vartheta^2}{2} V(x) \Big)   g''_\vartheta(x;0) \, \text{d} x\\
\! \!  &\! \! \! \!  & \! \!   + O(\vartheta^3) \\
& = & I_1(\vartheta) + \vartheta I_2(\vartheta) + \frac{\vartheta^2}{2} I_3(\vartheta) +O(\vartheta^3)
\end{eqnarray*}
(say). We have
\begin{eqnarray*}
I_1(\vartheta) & = & \frac{1}{(1+ 2 \delta)^{1/2}} + \vartheta \cdot \frac{\delta \sigma_1}{(1+ 2 \delta)^{3/2}} \\
& & + \vartheta^2 \bigg[ \frac{3 \delta^2\sigma_1^2}{2 (1+ 2 \delta)^{5/2}} + \frac{\delta (4\delta \mu_1^2 - 2 \sigma_1^2 + \sigma_2)}{2 (1+ 2 \delta)^{3/2}} - \frac{\delta \mu_1^2}{(1+ 2 \delta)^{1/2}}\bigg].
\end{eqnarray*}
Furthermore,
\begin{eqnarray*}
I_2(\vartheta) & = & \int \! {\rm e}^{-\delta x^2} g'_\vartheta(x;0) \, \text{d} x  + \delta \vartheta \int {\rm e}^{-\delta x^2} \big(\sigma_1^2x^2 + 2 \mu_1 x\big)  g'_\vartheta(x;0)  \, \text{d} x + O(\vartheta^2), \\
 I_3(\vartheta)& =  & \int \! {\rm e}^{-\delta x^2} g''_\vartheta(x;0) \, \text{d} x  + O(\vartheta).
\end{eqnarray*}
Recalling \eqref{defjk},
we thus obtain -- apart from a term which is $O(\vartheta^3)$ --
\begin{eqnarray}\nonumber
\! \! & \! \!  \! \!  & \! \! \int \exp\left(- \frac{\delta (x-\mu(\vartheta))^2}{\sigma^2(\vartheta)}\right) g(x;\vartheta) \, \text{d} x\\ \label{gleichungen2}
\! \! & \! \! = \! \! & \! \! \frac{1}{(1+ 2 \delta)^{1/2}} + \vartheta \bigg[ \frac{\delta \sigma_1}{(1+ 2 \delta)^{3/2}} + J_{1,0} \bigg]\\ \nonumber
\! \! & \! \! \! \! & \! \! + \vartheta^2 \bigg[ \frac{J_2}{2}+ \delta \sigma_1^2 J_{1,2} + 2 \delta \mu_1 J_{1,1} -\frac{\delta\left(\left(\delta+\frac12\right)(\sigma_2-2\mu_1^2)-\left(\frac{\delta}2+1\right)\sigma_1^2\right)}{(2\delta+1)^{5/2}}  \bigg].
\end{eqnarray}
Upon combining \eqref{erstesint} and \eqref{gleichungen2} and recalling \eqref{defdeltagam}, $b_{T_\beta}(\vartheta)$ figuring in
\eqref{integ1}, \eqref{integ2} takes the form
\[
b_{T_\beta}(\vartheta)  =  \Delta_\beta \vartheta^2 + O(\vartheta^3) \ \text{as} \ \vartheta \to 0,
\]
where
\begin{eqnarray*}
\Delta_\beta & = & D_0 + \frac{\beta^2}{(\beta^2+1)^{5/2}}\left(\left(\left(J_{1,0}-J_{1,2}\right)\sigma_1-2J_{1,1}\mu_1\right)\beta^2+J_{1,0}\sigma_1-2J_{1,1}\mu_1\right)\\ && +\frac{\beta^2}{(2\beta^2+1)^{5/2}}\left(\left(2\mu_1^2+\frac34\sigma_1^2\right)\beta^2+\mu_1^2\right),
\end{eqnarray*}
and $D_0$ and $J_{1,0}, J_{1,1}, J_{1,2}$ are defined in \eqref{defdj} and \eqref{defjk}, respectively.

\section{Approximations to solutions of the eigenvalue problem}\label{sec:approxev}
We now turn to the conditions \eqref{Bcondition1} and \eqref{Bcondition2}.
The limit null distribution of $T_{n,\beta}$, as $n \to \infty$, is given by the distribution of
\[ T_\beta :=  \int_{-\infty}^\infty Z^2(t) \, \varphi_\beta(t) \, {\rm d}t.
\]
Here, $Z$ is  a centred Gaussian random element of the Fr\'{e}chet space of continuous real-valued functions
having covariance kernel $K(s,t) = \mathbb{E}[Z(s)Z(t)]$, where
\begin{equation}\label{eq:KHW}
K(s,t)=\exp\left(-\frac{(s-t)^2}2\right)-\left(1+st+\frac{(st)^2}2\right)\exp\left(-\frac{s^2+t^2}2\right),\quad s,t\in\R
\end{equation}
(see Theorem 2.1 and Theorem 2.2 of \cite{HW:1997}). In fact, $Z$ may also be regarded as a Gaussian random element of the
separable Hilbert space $L^2$ (say) of (equivalence classes of) functions that are square integrable with respect to $\varphi_\beta(t) {\rm d}t$.
The distribution of $T_\beta$ is that of $\sum_{j=1}^\infty \lambda_j(\beta) N_j^2$, where $N_1, N_2, \ldots $ is a sequence of i.i.d. standard
normal random variables, and $\lambda_1(\beta), \lambda_2(\beta), \ldots $ is the sequence of positive eigenvalues of the integral operator
${\cal K}$  on $L^2$ defined by
\begin{equation*}
\mathcal{K}: L^2\rightarrow L^2,\quad f \mapsto \mathcal{K}f(s)=\int_{-\infty}^\infty K(s,t)f(t) \varphi_\beta(t)\,\mbox{d}t,\quad s\in\R.
\end{equation*}
Since $S_n$ figuring in \eqref{Bcondition1} equals $\sqrt{T_{n,\beta}}$, the function $G$ is the distribution function of $\widetilde{Z}:=\left(\sum_{j=1}^\infty \lambda_j(\beta) N_j^2\right)^{1/2}$.
From \cite{ZO:1961}, we thus  have
\[
\log(1-G(x)) = \log \PP(\widetilde{Z} >x) = \log \PP(\widetilde{Z}^2 > x^2) \sim - \frac{x^2}{2 \lambda_1(\beta)} \ \text{as} \  x \to \infty,
\]
where $\lambda_1(\beta)$ denotes the largest eigenvalue. Hence, the approximate Bahadur slope of the Epps--Pulley test statistic is given
by
\begin{equation}\label{appBEP}
c^*_{T_\beta}(\vartheta) = \frac{b_{T_\beta}(\vartheta)}{\lambda_1(\beta)}.
\end{equation}
Thus, one has to tackle  the so-called {\em eigenvalue problem}, i.e., to find positive values $\lambda$ and functions $f$ such that
$\mathcal{K}f =\lambda f$ or, in other words,  to solve the integral equation
\begin{equation}\label{intequ}
\int_{-\infty}^\infty K(s,t)f(t) \varphi_\beta(t)\,\mbox{d}t=\lambda f(s),\quad s\in\R.
\end{equation}
Since explicit solutions of such integral equations are only available in exceptional cases (for non-classical goodness-of-fit test statistics, see \cite{HN:2000} and
\cite {HN:2002}),  we employ a stochastic approximation method. This method is related to the quadrature method in the
classical literature on numerical mathematics (see \cite{B:1977}, chapter 3), and which can also be found in machine learning theory (see \cite{RW:2006}).
For the approximation of spectra of Hilbert Schmidt operators, see \cite{KG:2000}. To be specific, let $Y$ be a random variable having density $\varphi_\beta$.
Then (\ref{intequ}) reads
\begin{equation}\label{int:eq2}
\lambda f(s) = \E\big[K(s,Y)f(Y)\big],\quad s\in \R.
\end{equation}
An  empirical counterpart to (\ref{int:eq2}) emerges if we let $y_1,y_2,\ldots,y_N$, $N\in\N$, be independent realizations of $Y$.
An  approximation of the expected value in (\ref{int:eq2}) is then
\begin{equation}\label{int:eq3}
\E\big[K(s,Y)f(Y)\big]\approx \frac1N\sum_{j=1}^N K(s,y_j)f(y_j),\quad s\in \R.
\end{equation}
If we evaluate (\ref{int:eq3}) at the points $y_1,\ldots,y_n$, the result is
\begin{equation}\label{int:eq4}
\lambda f(y_i) = \frac1N\sum_{j=1}^N K(y_i,y_j)f(y_j),\quad i=1,\ldots,N,
\end{equation}
which is a system of $N$ linear equations. Writing $v=(f(y_1),\ldots,f(y_N))\in \R^N$ and $\widetilde{K}=(K(y_i,y_j)/N)_{i,j=1,\ldots,N}\in \R^{N\times N}$, we can rewrite (\ref{int:eq4})
according to
\begin{equation*}
\widetilde{K}v=\lambda v
\end{equation*}
in matrix form, from which the (approximated) eigenvalues $\lambda_1,\ldots,\lambda_N$ can be computed explicitly. Note that for each
eigenvalue $\lambda_j$ we have an eigenvector  $v_j\in\R^N$ (say), the components of which are the (approximated) values of the eigenfunctions (say) $f_j$ computed at $y_1,\ldots,y_N$.

The simulation of eigenvalues was performed in the  statistical computing language \texttt{R}, see \cite{R:2021}. As parameters for the
simulation we chose $N=1000$, and we considered the tuning parameters $\beta\in\{0.25,0.5,0.75,1,2,3,5,10\}$. Each entry in Table 1 stands for the mean of 10 simulation runs.
\begin{table}[t]
\centering
\begin{tabular}{r|rrrrrrrr}
 $\lambda\backslash \beta$ & 0.25 & 0.5 & 0.75 & 1 & 2 & 3 & 5 & 10 \\
  \hline
$\lambda_1$ & 0.00040 & 0.01065 & 0.03829 & 0.07507 & 0.15207 & 0.16149 & 0.13552 & 0.08791 \\
  $\lambda_2$ & 0.00003 & 0.00304 & 0.01735 & 0.04454 & 0.12921 & 0.14577 & 0.12606 & 0.08178 \\
  $\lambda_3$ & 0.00000 & 0.00021 & 0.00220 & 0.00846 & 0.04894 & 0.07676 & 0.08703 & 0.06879 \\
  $\lambda_4$ & 0.00000 & 0.00004 & 0.00076 & 0.00417 & 0.03966 & 0.06642 & 0.07997 & 0.06459 \\
  $\lambda_5$ & 0.00000 & 0.00000 & 0.00011 & 0.00098 & 0.01692 & 0.03755 & 0.05678 & 0.05518 \\
\end{tabular}
\caption{Approximate first five eigenvalues of $\mathcal{K}$ for different weight functions $\varphi_\beta$, each entry is the mean of 10 simulation runs}\label{tab:ev}
\end{table}

\section{Alternatives}\label{sec:alternatives}
As in Milo\u{s}evi\'{c} et al. (\cite{MNO:2021}), we consider the following close alternatives:
\begin{itemize}
\item a Lehmann alternative with density
\[
g_1(x;\vartheta) = (1+\vartheta)\Phi^\vartheta(x)\varphi(x);
\]
\item a first Ley-Paindaveine alternative with density (see \cite{LP:2009})
\[
g_2(x;\vartheta) = \varphi(x) {\rm e}^{-\vartheta(1-\Phi(x))}(1+\vartheta \Phi(x));
\]
\item a second Ley-Paindaveine alternative with density (see \cite{LP:2009})
\[ g_3(x;\vartheta) = \varphi(x)(1-\vartheta \pi \cos(\pi \Phi(x)));
\]
\item a contamination alternative (with $\text{N}(\mu,\sigma^2)$ for several pairs $(\mu,\sigma^2) \neq (0,1)$) with density
\[
g_4^{[\mu,\sigma^2]}(x;\vartheta) = (1- \vartheta)\varphi(x) + \frac{\vartheta}{\sigma}\varphi\left(\frac{x-\mu}{\sigma} \right).
\]
\end{itemize}
Like in Milo\u{s}evi\'{c} et al. (\cite{MNO:2021}), we computed the local (as $\vartheta \to 0$) relative approximate Bahadur efficiencies with respect to
the likelihood ratio test (LRT). The LRT is the best test regarding exact Bahadur effiency, and it is often used as a benchmark test.
Table 2 displays the local approximate Bahadur efficiencies of $T_{n,\beta}$ with respect to the LRT, for each of the six alternatives
considered in \cite{MNO:2021}, and for $\beta \in \{0.25,0.5, 0.75, 1,2,4,5,10\}$. A comparison with Table 1 of \cite{MNO:2021} shows that the Epps--Pulley test
with $\beta = 0.5$ dominates the Kolmogorov--Smirnov test for each of the six alternatives, and for $\beta = 1$, $\beta = 2$ and $\beta =3$, it outperforms
the tests of Cram\'{e}r--von Mises, the Watson variation of this test, and the Watson--Darling variation of the Kolmogorov--Smirnov test, respectively.
If $\beta = 0.75$, the Epps--Pulley test dominates the Anderson--Darling test for each of the alternatives with the exception of the final contamination alternative.
As a conclusion, the test of Epps and Pulley should receive more attention as a test for normality.

\begin{table}[t]
\centering
\begin{tabular}{r|rrrrrrrr}
 Alt.$\backslash \beta$ & 0.25 & 0.5 & 0.75 & 1 & 2 & 3 & 5 & 10 \\
  \hline
Lehmann & 0.996 & 0.895  & 0.854 & 0.743 & 0.514 & 0.406 & 0.328 & 0.267 \\
 1st Ley-Paindaveine & 0.947 & 0.944 & 0.998 & 0.937 & 0.745 & 0.612 & 0.507 & 0.417 \\
 2nd Ley-Paindaveine & 0.824 & 0.872 & 0.986 & 0.981 & 0.881 & 0.754 & 0.641 & 0.533 \\
  Contamination with N(1,1) & 0.760 & 0.649 & 0.592 & 0.499 & 0.328 & 0.255 & 0.205 & 0.166 \\
  Contamination with N(0.5,1) & 0.945 & 0.824 & 0.766 & 0.654 & 0.438 & 0.343 & 0.276 & 0.224 \\
  Contamination with N(0,0.5) & 0.084 & 0.267 & 0.474 & 0.587 & 0.675 & 0.606 & 0.526 & 0.442 \\
\end{tabular}
\caption{Approximate local Bahadur efficiency of $T_{n,\beta}$ with respect to the LRT }\label{tab:ev}
\end{table}

\bibliography{lit_BHEP_EV}   
\bibliographystyle{abbrv}      
\vspace{5mm}
B. Ebner and N. Henze, \\
Institute of Stochastics, \\
Karlsruhe Institute of Technology (KIT), \\
Englerstr. 2, D-76133 Karlsruhe. \\
E-mail: {\texttt Bruno.Ebner@kit.edu}\\
E-mail: {\texttt Norbert.Henze@kit.edu}

\end{document}